\documentclass[a4paper,10pt, reqno]{amsart}
\usepackage[backref]{hyperref}
\usepackage{a4wide}
\usepackage{amsthm}
\usepackage{amssymb}
 \newtheorem{theorem}{Theorem}
 
 \newtheorem{lemma}{Lemma}

\begin{document}
\title{The square sieve and the large sieve with square moduli}

\author{Stephan Baier}
\address{Stephan Baier, School of Mathematics, Tata Institute of Fundamental Research, 1 Dr. Homi Bhaba Road, Colaba, Mumbai 400005, India}

\email{sbaier@math.tifr.res.in}


\maketitle

{\bf Abstract:} We give a short alternative proof using Heath-Brown's square sieve of a bound of the author for the large sieve with square moduli.  

\section{Main result} 
The large sieve with square moduli and, more generally, power moduli was studied by L. Zhao and the author in a number of papers, both independently and in joint work. The best result for square moduli obtained so far is \cite[Theorem 1]{BZ} which asserts the following.

\begin{theorem} \label{bazo} Let $\varepsilon>0$. Then for any $M\in \mathbb{Z}$, $N\in \mathbb{N}$, $Q\ge 1$ and sequence of complex numbers $\left(a_n\right)_{n\in\mathbb{Z}}$,  we have
\begin{equation} \label{what}
\sum\limits_{q\le Q} \sum\limits_{\substack{a=1\\ (a,q)=1}}^{q^2} \left|S\left(\frac{a}{q^2}\right)\right|^2= O\left((NQ)^{\varepsilon}\left(Q^3+N+\min\left\{N\sqrt{Q},\sqrt{N}Q^2\right\}\right)Z\right),
\end{equation}
where 
\begin{equation} \label{nots}
S(\alpha):=\sum\limits_{n=M+1}^{M+N} a_ne(n\alpha) \quad \mbox{and} \quad Z:=\sum\limits_{n=M+1}^{M+N} \left|a_n\right|^2.
\end{equation}
\end{theorem}

In \cite{Zh}, Zhao proved \eqref{what} with $N\sqrt{Q}+\sqrt{N}Q^2$ in place of the minimum of these terms using Fourier analysis. By combinatorial considerations, the author of the present paper then showed in \cite{Ba} that the term $N\sqrt{Q}$ in the above sum can be omitted. Finally, combining their methods and making further refinements, both authors together succeeded in proving \eqref{what}. We note that the two terms in the minimum coincide and give a contribution of $Q^{7/2}$ if $Q^3=N$. It was conjectured by Zhao that \eqref{what} should hold without the minimum term. 

The purpose of this paper is to give a short alternative proof of the bound
\begin{equation} \label{goal}
\sum\limits_{q\le Q} \sum\limits_{\substack{a=1\\ (a,q)=1}}^{q^2} \left|S\left(\frac{a}{q^2}\right)\right|^2 = O\left( (NQ)^{\varepsilon}\left(Q^3+N+\sqrt{N}Q^2\right)Z\right),
\end{equation}
previously obtained in \cite{Ba},
using Heath-Brown's square sieve. We note that the left-hand side of \eqref{goal} is bounded by 
\begin{equation} \label{majo1}
\sum\limits_{q\le Q} \sum\limits_{\substack{a=1\\ (a,q)=1}}^{q^2} \left|S\left(\frac{a}{q^2}\right)\right|^2 \ll \sum\limits_{q\le Q^2} \sum\limits_{\substack{a=1\\ (a,q)=1}}^{q} \left|S\left(\frac{a}{q}\right)\right|^2 \ll \left(Q^4+N\right)Z
\end{equation}
using the classical large sieve inequality (see, for example, Zhao's thesis \cite{Zh} for reference), and also by 
\begin{equation} \label{majo2}
\sum\limits_{q\le Q} \sum\limits_{\substack{a=1\\ (a,q)=1}}^{q^2} \left|S\left(\frac{a}{q^2}\right)\right|^2 \ll Q(N+Q^2)Z,
\end{equation}
which is obtained by summing over $q$ the large sieve bound 
$$
\sum\limits_{\substack{a=1}}^{q} \left|S\left(\frac{a}{q}\right)\right|^2 \ll \left(q+N\right)Z
$$
for individual moduli (which can be seen by just expanding the square and re-arranging summations). Since Zhao's conjecture follows from \eqref{majo1} if $Q\le N^{1/4+\varepsilon}$ and from \eqref{majo2} if $Q\ge N^{1/2-\varepsilon}$, we shall assume that
\begin{equation} \label{Qcond}
N^{1/4+\varepsilon}<Q<N^{1/2-\varepsilon}
\end{equation}
throughout this paper. We note that \eqref{goal} is the same as \eqref{what} if $Q\le N^{1/3}$. \\
$ $\\
{\bf Notation and conventions:} \medskip\\
(i) We keep the notations in Theorem \ref{bazo}, in particular those in \eqref{nots}, throughout this paper. \\
(ii) $n=\Box$ means that $n$ is a square.\\ 
(iii) $d(n):=\sum_{m|n} 1$ denotes the divisor function.\\
(iv) $\varepsilon$ is an arbitrarily small positive number which can change from line to line.\\ 

{\bf Acknowledgement.} The author wishes to thank the Tata Institute of Fundamental Research in Mumbai (India) for its warm hospitality, excellent working conditions and financial support by an ISF-UGC grant.

\section{Preliminaries}
Similarly as in \cite{Ba}, our starting point will be the following Lemma which is a direct consequence of \cite[Theorem 2.1]{Mo}. 

\begin{lemma} \label{largesieve}
Assume that $Q\ge 1$, $N\ge 1$ and $0<\Delta\le 1$. Then
\begin{equation} \label{plugging}
\begin{split}
& \sum\limits_{Q<q\le 2Q} \sum\limits_{\substack{a=1\\ (a,q)=1}}^{q^2} \left|S\left(\frac{a}{q^2}\right)\right|^2 \ll \left(N+\Delta^{-1}\right)Z \cdot
\max\limits_{\alpha\in \mathbb{R}}  P\left(\alpha,\Delta\right),
\end{split}
\end{equation}
where 
\begin{equation} \label{Pi}
P(\alpha,\Delta):=\mathop{\sum\limits_{Q<q\le 2Q} \sum\limits_{\substack{a=1\\ (a,q)=1}}^{q^2}}_{|a/q^2-\alpha|\le \Delta} 1.
\end{equation}
\end{lemma}

To detect squares, we shall employ Heath-Brown's square sieve (see \cite[Theorem 1]{HB}). 

\begin{lemma} \label{squaresieve}
Let $\mathcal{P}$ be a set of $P$ primes. Suppose that $w : \mathbb{Z} \rightarrow \mathbb{R}$ is a function satisfying $w(n)\ge 0$ for all $n\in \mathbb{N}$ and $w(n)=0$ if $n=0$ or $|n|\ge e^P$.  Then
\begin{equation*}
\sum\limits_{n=1}^{\infty} w\left(n^2\right) \ll P^{-1}\sum\limits_{n\in \mathbb{Z}} w(n) +P^{-2} \sum\limits_{\substack{p_1,p_2\in \mathcal{P}\\ p_1\not=p_2}} \left|\sum\limits_{n\in \mathbb{Z}} w(n) \left(\frac{n}{p_1p_2}\right)\right|,
\end{equation*}
where $\left(\frac{n}{p_1p_2}\right)$ is the Jacobi symbol.
\end{lemma}

\section{Division into major and minor arcs} 
Dividing the $q$-range in \eqref{what} into dyadic intervals, it suffices to prove that
\begin{equation} \label{prove}
\sum\limits_{Q<q\le 2Q} \sum\limits_{\substack{a=1\\ (a,q)=1}}^{q^2} \left|S\left(\frac{a}{q^2}\right)\right|^2= O\left((NQ)^{\varepsilon}\left(Q^3+N+\sqrt{N}Q^2\right)Z\right). 
\end{equation}
Our task is now to estimate the term $P(\alpha,\Delta)$ in Lemma \ref{largesieve} for any given $\alpha\in \mathbb{R}$. We shall choose 
\begin{equation} \label{choose}
\Delta:=\frac{1}{N}.
\end{equation}
First, we first consider $\alpha$ in a set of major arcs, defined by
\begin{equation*}
\mathfrak{M}:= \bigcup\limits_{v\le 1/(500Q^2\Delta)}\bigcup\limits_{\substack{u=1\\ (u,v)=1}}^v \left[\frac{u}{v}-\frac{1}{10Q^2v},\frac{u}{v}+\frac{1}{10Q^2v}\right].  
\end{equation*}
If $\alpha\in \mathfrak{M}$, then there exist $u,v\in \mathbb{Z}$ such that $1\le v\le 1/(500Q^2\Delta)$, $(u,v)=1$ and $|u/v-\alpha|\le 1/(10Q^2v)$ and hence
\begin{equation} \label{Mcase}
P\left(\alpha,\Delta\right)\le P\left(\frac{u}{v},\frac{1}{5Q^2v}\right) = \mathop{\sum\limits_{Q<q\le 2Q} \sum\limits_{\substack{a=1\\ (a,q)=1}}^{q^2}}_{|a/q^2-u/v|\le 1/(5Q^2v)} 1 \le
\mathop{\sum\limits_{Q<q\le 2Q} \sum\limits_{\substack{a=1\\ (a,q)=1}}^{q^2}}_{|av-uq^2|<1} 1 = 
\mathop{\sum\limits_{Q<q\le 2Q} \sum\limits_{\substack{a=1\\ (a,q)=1}}^{q^2}}_{a/q^2=u/v} 1 = 0
\end{equation}
since $(a,q^2)=1=(u,v)$ and $v\le 1/(500Q^2\Delta)<Q^2< q^2$ by \eqref{Qcond} if $N$ is large enough. Hence, the major arcs don't contribute. 

In the remainder, we consider the case when $\alpha$ is in the set of minor arcs, defined by 
$$
\mathfrak{m}:= \mathbb{R}\setminus \mathfrak{M}.
$$

\section{Application of the square sieve}
By Dirichlet's approximation theorem, there exist integers $b$ and $r$ such that 
\begin{equation} \label{brcondis}
1\le r\le 500 Q^2, \quad (b,r)=1 \quad \mbox{and} \quad \left| \frac{b}{r} -\alpha\right| \le \frac{1}{500Q^2r}.
\end{equation} 
If $r\le 1/(500Q^2\Delta)$, then it follows that $\alpha\in \mathcal{M}$. Thus, \eqref{brcondis} can be replaced by 
\begin{equation} \label{brcondisnew}
1/(500Q^2\Delta)<r\le 500 Q^2, \quad (b,r)=1, \quad \mbox{and} \quad \left| \frac{b}{r} -\alpha\right| \le \Delta.
\end{equation}
It follows that 
\begin{equation} \label{firstbound}
P(\alpha,\Delta)\le P\left(\frac{b}{r},2\Delta\right).
\end{equation}

Let $\Phi_1$ and $\Phi_2$ be infinitely differentiable compactly supported functions from $\mathbb{R}$ to $\mathbb{R}^+$, supported in the intervals $[1/2, 5]$ and $[-10, 10]$ and bounded below by 1 on the intervals $[1, 4]$ and $[-4, 4]$, respectively. Then
\begin{equation} \label{Pi}
P\left(\frac{b}{r},2\Delta\right)\ll \sum\limits_{q\in \mathbb{Z}} \Phi_1\left(\frac{q^2}{Q^2}\right)\cdot \sum\limits_{a\in \mathbb{Z}} \Phi_2\left(\frac{a-q^2b/r}{Q^2 \Delta}\right).
\end{equation}
Let 
\begin{equation} \label{parameter}
R> (QN)^{\varepsilon}
\end{equation} 
be a parameter, to be fixed later, and 
\begin{equation} \label{mathP}
\mathcal{P}:= \left\{ p\in \mathbb{P} \ :\ R<p\le 2R \mbox{ and } p\nmid r\right\}, 
\end{equation}
where $\mathbb{P}$ is the set of all primes. In the notation of Lemma \ref{squaresieve}, we have
\begin{equation} \label{P}
P:=\sharp \mathcal{P}=\pi(2R)-\pi(R)-\omega(r) \sim \frac{R}{\log R}.
\end{equation} 
Now applying the square sieve, Lemma \ref{squaresieve}, to the right-hand side of \eqref{Pi}, we get
\begin{equation} \label{aftersquaresieve}
\begin{split}
P\left(\frac{b}{r},2\Delta\right) \ll & \frac{1}{P} \cdot \sum\limits_{n\in \mathbb{Z}} \Phi_1\left(\frac{n}{Q^2}\right) \cdot  \sum\limits_{a\in \mathbb{Z}} \Phi_2\left(\frac{a-nb/r}{Q^2 \Delta}\right) + \\ &
\frac{1}{P^{2}} \cdot \sum\limits_{\substack{p_1,p_2\in \mathcal{P}\\ p_1\not=p_2}} \left| \sum\limits_{n\in \mathbb{Z}} \Phi_1\left(\frac{n}{Q^2}\right) \cdot \left(\frac{n}{p_1p_2}\right) \cdot \sum\limits_{a\in \mathbb{Z}} \Phi_2\left(\frac{a-bn/r}{Q^2\Delta}\right)  \right|.
\end{split}
\end{equation}

The first double sum over $n$ and $a$ on the right-hand side of \eqref{aftersquaresieve} can be estimated by
\begin{equation*} 
\begin{split}
& \sum\limits_{n\in \mathbb{Z}} \Phi_1\left(\frac{n}{Q^2}\right)  \cdot \sum\limits_{a\in \mathbb{Z}} \Phi_2\left(\frac{a-nb/r}{Q^2 \Delta}\right) \le
\mathop{\sum\limits_{Q^2/2\le n\le 5Q^2} \ \sum\limits_{a\in\mathbb{Z}}}_{|a/n-b/r|\le 20\Delta} \ 1 \\ = & 
 \mathop{\sum\limits_{Q^2/2\le n\le 5Q^2} \ \sum\limits_{a\in \mathbb{Z}}}_{\substack{(a,n)\le 2500Q^4\Delta\\ |a/n-b/r|\le 20\Delta}} \ 1+
\mathop{\sum\limits_{Q^2/2\le n\le 5Q^2} \ \sum\limits_{a\in \mathbb{Z}}}_{\substack{(a,n)> 2500Q^4\Delta\\ |a/n-b/r|\le 20\Delta}} \ 1\\ \le
& \sum\limits_{d\le 2500Q^4\Delta} \ \mathop{\sum\limits_{Q^2/(2d)\le n_1\le 5Q^2/d}\ \sum\limits_{a_1\in \mathbb{Z}}}_{\substack{(a_1,n_1)=1\\ |a_1/n_1-b/r|\le 20\Delta}}\  1+
\mathop{\sum\limits_{n_1\le 1/(500Q^2\Delta)} \ \sum\limits_{a_1\in \mathbb{Z}}}_{\substack{(a_1,n_1)=1\\ |a_1/n_1-b/r|\le 20\Delta}} \ \sum\limits_{Q^2/(2n_1)\le d\le 5Q^2/n_1}\ 1.
\end{split}
\end{equation*}
Since $|a_1/n_1-a_2/n_2|\ge d^2/(25Q^4)$ whenever $Q^2/(2d)\le n_1,n_2\le 5Q^2/d$, $(a_1,n_1)=1=(a_2,n_2)$ and $a_1/n_1\not=a_2/n_2$, it follows that
\begin{equation*} 
 \sum\limits_{d\le 2500Q^4\Delta} \ \mathop{\sum\limits_{Q^2/(2d)\le n_1\le 5Q^2/d}\ \sum\limits_{a_1\in \mathbb{Z}}}_{\substack{(a_1,n_1)=1\\ |a_1/n_1-b/r|\le 20\Delta}}\  1\ll 
 \sum\limits_{d\le 2500Q^4\Delta}\left(1+\frac{Q^4\Delta}{d^2}\right) \ll Q^4\Delta.
\end{equation*}
Further, 
\begin{equation*} 
\mathop{\sum\limits_{n_1\le 1/(500Q^2\Delta)} \ \sum\limits_{a_1\in \mathbb{Z}}}_{\substack{(a_1,n_1)=1\\ |a_1/n_1-b/r|\le 20\Delta}} \ \sum\limits_{Q^2/(2n_1)\le d\le 5Q^2/n_1}\ 1 =0
\end{equation*}
since $|a_1/n_1-b/r|\le 20\Delta$ implies $|a_1/n_1-\alpha|\le 21\Delta$ by \eqref{brcondisnew}, and hence, $\alpha\in \mathfrak{M}$, which is not the case. Altogether, we thus obtain
\begin{equation} \label{firstsum} 
\sum\limits_{n\in \mathbb{Z}} \Phi_1\left(\frac{n}{Q^2}\right)  \cdot \sum\limits_{a\in \mathbb{Z}} \Phi_2\left(\frac{a-nb/r}{Q^2 \Delta}\right) \ll Q^4\Delta.
\end{equation}

\section{Application of Poisson summation}
Now we estimate the second double sum over $n$ and $a$ on the right-hand side of \eqref{aftersquaresieve}. We first split the summation over $n$ into residue classes modulo $p_1p_2r$, getting
\begin{equation} \label{p1p2rsplit}
\begin{split}
& \sum\limits_{n\in \mathbb{Z}} \Phi_1\left(\frac{n}{Q^2}\right) \cdot \left(\frac{n}{p_1p_2}\right) \cdot \sum\limits_{a\in \mathbb{Z}} \Phi_2\left(\frac{a-bn/r}{Q^2\Delta}\right)\\
= & \sum\limits_{f=1}^{p_1p_2r} \left(\frac{f}{p_1p_2}\right) \cdot  \sum\limits_{(x_1,x_2)\in \mathbb{Z}^2} \Phi_1\left(\alpha_1(f+mx_1)\right) \cdot \Phi_2\left(\beta_2x_2+\alpha_2(f+mx_1)\right),
\end{split}
\end{equation}
where 
\begin{equation} \label{varia}
\alpha_1:=\frac{1}{Q^2}, \quad \beta_2:=\frac{1}{Q^2\Delta}, \quad \alpha_2:=-\frac{b}{rQ^2\Delta}, \quad m:=p_1p_2r.
\end{equation}
Now we apply the two-dimensional Poisson summation 
\begin{equation} \label{poisson}
\sum\limits_{(x_1,x_2)\in \mathbb{Z}} \Phi(x_1,x_2) = \sum\limits_{(t_1,t_2)\in \mathbb{Z}} \hat\Phi(t_1,t_2)
\end{equation}
to the function
\begin{equation} \label{Phide}
\Phi(x_1,x_2):=\Phi_1\left(\alpha_1(f+mx_1)\right) \cdot \Phi_2\left(\beta_2x_2+\alpha_2(f+mx_1)\right)
\end{equation}
whose Fourier transform is 
\begin{equation} \label{ftransform}
\hat\Phi(t_1,t_2) = \frac{1}{\beta_2\alpha_1m} \cdot e\left(\frac{lt_1}{m}\right) \cdot \hat\Phi_1\left(\frac{t_1}{\alpha_1m}-\frac{t_2\alpha_2}{\alpha_1\beta_2}\right) \cdot \hat\Phi_2\left(\frac{t_2}{\beta_2}\right).
\end{equation}
Combining \eqref{p1p2rsplit}, \eqref{varia}, \eqref{poisson}, \eqref{Phide} and \eqref{ftransform}, we deduce that 
\begin{equation} \label{endlich} 
\begin{split}
& \sum\limits_{n\in \mathbb{Z}} \Phi_1\left(\frac{n}{Q^2}\right) \cdot \left(\frac{n}{p_1p_2}\right) \cdot \sum\limits_{a\in \mathbb{Z}} \Phi_2\left(\frac{a-bn/r}{Q^2\Delta}\right)\\ \\
= & \frac{Q^4 \Delta}{p_1p_2r} \cdot \sum\limits_{c\in \mathbb{Z}}\hat\Phi_2\left(cQ^2 \Delta\right)\cdot   \sum\limits_{s\in \mathbb{Z}} \hat\Phi_1\left(\frac{sQ^2}{p_1p_2r}\right)
\sum\limits_{f=1}^{p_1p_2r} \left(\frac{f}{p_1p_2}\right) \cdot e\left(\frac{(s-cbp_1p_2)f}{p_1p_2r}\right).
\end{split}
\end{equation}
Since $(r,p_1p_2)=1$, we have  
\begin{equation} \label{dieexpsum}
\begin{split}
\sum\limits_{f=1}^{p_1p_2r} \left(\frac{f}{p_1p_2}\right) \cdot e\left(\frac{(s-cbp_1p_2)f}{p_1p_2r}\right)= & \sum\limits_{g=1}^{p_1p_2} \sum\limits_{h=1}^r 
\left(\frac{gr}{p_1p_2}\right) \cdot e\left(\frac{(s-cbp_1p_2)(gr+hp_1p_2)}{p_1p_2r}\right)\\
= & \sum\limits_{g=1}^{p_1p_2} \left(\frac{gr}{p_1p_2}\right) \cdot e\left(\frac{sg}{p_1p_2}\right) \cdot \sum\limits_{h=1}^r e\left(\frac{(s-cbp_1p_2)h}{r}\right)\\
= &  r\tau_{p_1p_2}\cdot \left(\frac{rs}{p_1p_2}\right) \cdot \begin{cases} 1 & \mbox{ if } cp_1p_2\equiv s\overline{b} \bmod{r} \\ 0 & \mbox{ otherwise,}\end{cases}
\end{split}
\end{equation}
where $\tau_{p_1p_2}$ is the Gauss sum for the Jacobi symbol modulo $p_1p_2$.  
Combining \eqref{endlich} and \eqref{dieexpsum}, and using the triangle inequality and the well-known equation $\left|\tau_{p_1p_2}\right|=\sqrt{p_1p_2}$, we get
\begin{equation} \label{end} 
\begin{split}
& \left| \sum\limits_{n\in \mathbb{Z}} \Phi_1\left(\frac{n}{Q^2}\right) \cdot \left(\frac{n}{p_1p_2}\right)\cdot  \sum\limits_{a\in \mathbb{Z}} \Phi_2\left(\frac{a-bn/r}{Q^2\Delta}\right)\right| \\
\le & \frac{Q^4 \Delta}{\sqrt{p_1p_2}} \cdot \mathop{\sum\limits_{c\in \mathbb{Z}} \sum\limits_{s\in \mathbb{Z}}}_{cp_1p_2 \equiv s\overline{b} \bmod{r}} 
\left| \hat\Phi_1\left(\frac{sQ^2}{p_1p_2r}\right) \cdot \hat\Phi_2\left(cQ^2 \Delta\right)\right|.
\end{split}
\end{equation} 

\section{Completion of the proof}
Let $T=(QRN)^{\varepsilon}$. Summing the right-hand side of \eqref{end} over $p_1,p_2$, and using the rapid decays of $\hat\Phi_1$ and $\hat\Phi_2$, we get
\begin{equation}  \label{above}
\begin{split}
& \sum\limits_{\substack{p_1,p_2\in \mathcal{P}\\ p_1\not= p_2}} \frac{Q^4 \Delta}{\sqrt{p_1p_2}} \cdot \mathop{\sum\limits_{c\in \mathbb{Z}} \sum\limits_{s\in \mathbb{Z}}}_{cp_1p_2 \equiv s\overline{b} \bmod{r}} 
\left| \hat\Phi_1\left(\frac{sQ^2}{p_1p_2r}\right) \cdot \hat\Phi_2\left(cQ^2 \Delta\right)\right|\\
\ll &  1+\frac{Q^4 \Delta}{R}\cdot  \mathop{\sum\limits_{R^2<m\le (2R)^2} \sum\limits_{|c|\le T/(Q^2\Delta)} \sum\limits_{|s|\le TR^2r/Q^2}}_{cm \equiv s\overline{b} \bmod{r}} 1\\
\ll & 1+\frac{Q^4 \Delta}{R} \cdot \left(R^2\cdot \left(1+\frac{TR^2}{Q^2}\right)+\mathop{\sum\limits_{|s|\le TR^2r/Q^2} \sum\limits_{R^2<m\le (2R)^2} \sum\limits_{1\le |c|\le T/(Q^2\Delta)}}_{cm \equiv s\overline{b} \bmod{r}}  1\right)\\
\ll & 1+Q^4 \Delta R + TQ^2 \Delta R^3+\frac{Q^4 \Delta}{R}\cdot \sum\limits_{|s|\le TR^2r/Q^2} \sum\limits_{\substack{1\le t\le 4R^2T/(Q^2\Delta)\\ t\equiv s\overline{b}\bmod{r}}} d(t)\\
\ll & 1+Q^4 \Delta R + TQ^2 \Delta R^3+\frac{TQ^4 \Delta}{R}\cdot \left(1+\frac{TR^2r}{Q^2}\right) \cdot \left(1+ \frac{R^2T}{rQ^2\Delta}\right) \\
\ll &  1+T^3\left(Q^4 \Delta R + Q^2 \Delta R^3+\frac{Q^4 \Delta}{R} + Q^2\Delta Rr + \frac{Q^2R}{r}+ R^3\right) \\
\ll & 1+T^3\left(Q^4 \Delta R + Q^2 \Delta R^3+R^3\right),
\end{split}
\end{equation} 
where for the last line, we use the bounds for $r$ in \eqref{brcondisnew} and \eqref{parameter}. Combining \eqref{P}, \eqref{aftersquaresieve}, \eqref{firstsum}, \eqref{end} and \eqref{above}, we get
\begin{equation*} 
\begin{split}
P\left(\frac{b}{r},2\Delta\right) \ll  (QRN)^{\varepsilon}\cdot \left(\frac{Q^4\Delta}{R} + Q^2 \Delta R + R\right).
\end{split}
\end{equation*}
Now we choose
$$
R:=Q^2\sqrt{\Delta} 
$$
which is consistent with \eqref{parameter} by \eqref{Qcond} and \eqref{choose} and gives
$$
P\left(\frac{b}{r},2\Delta\right) \ll  (QN)^{\varepsilon}\cdot \left(Q^2\Delta^{1/2} + Q^4 \Delta^{3/2}\right).
$$
Plugging this into \eqref{plugging} and using \eqref{Qcond} and \eqref{choose} proves \eqref{goal}. $\Box$

\end{document}